\documentclass[final,3p,11pt]{elsarticle}
\setcitestyle{citesep={;}}
\usepackage{graphicx}
\usepackage{bm}
\usepackage{timet}
\usepackage{setspace}
\usepackage{soul} 

\usepackage{xcolor} 
\usepackage{xfrac}
\usepackage{lineno}
\usepackage[colorlinks]{hyperref}
\hypersetup{
     colorlinks   = true,
     citecolor    = blue
}
\usepackage[section]{placeins}
\usepackage{algorithmic}
\usepackage{algorithm}
\usepackage{amssymb, amsmath}
\journal {Computer \& Mathematics with Applications}
\begin{document}
\begin{frontmatter}
\title{\begin{myfont} \LARGE \textbf{A stabilized radial basis-finite difference (RBF-FD) method with hybrid kernels} \end{myfont}} 
\author[label1]{\begin{myfont} Pankaj K Mishra\corref{cor1}\end{myfont}}
\address[label1]{\begin{myfont}Department of Mathematics, Hong Kong Baptist University, Kowloon Tong, Hong Kong\end{myfont}}
\ead{pankajkmishra01@gmail.com}
\cortext[cor1]{\begin{myfont}Corresponding Author\end{myfont}}
\ead[url]{https://mishrapk.com}
\author[label2]{\begin{myfont}Gregory E Fasshauer\end{myfont}}
\address[label2]{\begin{myfont}Department of Applied Mathematics \& Statistics, Colorado School of Mines, Golden, USA\end{myfont}}

\author[label3]{\begin{myfont}Mrinal K Sen\end{myfont}}
\address[label3]{\begin{myfont}Institute for Geophysics, University of Texas at Austin, USA\end{myfont}}
\author[label1]{\begin{myfont}Leevan Ling\end{myfont}}
\end{frontmatter}
\begin{myfont}

\section*{Abstract}
\noindent Recent developments have made it possible to overcome grid-based limitations of  finite difference (FD) methods by adopting the kernel-based meshless framework using radial basis functions (RBFs). Such an approach provides a meshless implementation and is referred to as the radial basis-generated finite difference (RBF-FD) method. In this paper, we propose a stabilized RBF-FD approach with a hybrid kernel, generated through a hybridization of the Gaussian and cubic RBF. This hybrid kernel was found to improve the condition of the system matrix, consequently, the linear system can be solved with direct solvers which leads to a significant reduction in the computational cost as compared to standard RBF-FD methods coupled with present stable algorithms. Unlike other RBF-FD approaches, the eigenvalue spectra of differentiation matrices were found to be stable irrespective of irregularity, and the size of the stencils. As an application, we solve the frequency-domain acoustic wave equation in a 2D half-space. In order to suppress spurious reflections from truncated computational boundaries, absorbing boundary conditions have been effectively implemented.
\end{myfont}
\begin{myfont}
\onecolumn 
\section{Introduction}
Numerical modeling is an indispensable step in understanding complex processes in science and engineering, occurring either at significantly large scale or micro-scale, where an interpretation with direct measurements is not always possible. Milestone developments in such numerical approaches are: finite difference (FD), finite element (FE), finite volume (FV), and pseudospectral (PS) methods. Each of these methods approximates the solution of the governing equations on a distribution of nodes or elements, which are fixed through a grid of points (or \textit{mesh}). In the last two decades, however, a new family of \textit{meshless} numerical methods has drawn the attention of numerical modelers. These methods do not require physical connections between the nodes but rather use the interaction of each node with all or a certain number of neighbor nodes through a specified kernel.

Kernel-based finite difference methods represent a local meshless approach which has been shown to work well for large scale problems. Radial basis functions have been popular choices for the kernels of such methods. Therefore one frequently refers to this approach as the radial basis-finite difference method (RBF-FD). Application of RBFs to compute derivatives on unstructured grids was first introduced by \cite{Tolstykh2000} and then formally proposed as RBF-FD approach in \cite{Tolstykh2003,Wright2003}. Since then, this method has been continuously improved \cite{Wright2006,Bayona2010,Fornberg2011,Fornberg2015} and applied to numerical modeling for various processes including convection-diffusion \cite{Chandhini2007,Stevens2009}, Navier-Stokes \cite{Chinchapatnam2009}, atmospheric global electric circuit \cite{Bayona2010}, shallow water simulation \cite{Flyer2012}, reaction-diffusion on surfaces \cite{Shankar2015,Shankar2017}, and time-domain elastic wave propagation in 2D isotropic media \cite{Brad2015} and heat flow \cite{Martin2017}.

Use of infinitely smooth radial kernels, like Gaussians, theoretically provides a spectrally convergent meshless method. The Gaussian RBF interpolant converges to a polynomial interpolant when the shape parameter tends to zero \cite{Driscoll2002413,Fornberg200437,Fassbook2007}. Taking into consideration this connection between polynomials and RBFs, \cite{Larsson2005} have shown that it is possible to achieve higher accuracy in meshless approximation by using nearly flat Gaussian kernels. In modern conventions, the shape parameter is inversely proportional to the average distance between the node points. Use of larger shape parameters may result in an approximation that is similar to overfitting. Therefore, it is recommended to keep the shape parameter on the smaller side. A `small' \footnote{$\epsilon$ in the Table \ref{tab:rbf}} shape parameter thus becomes typical for kernel-based meshless algorithms to offer a little bit extra smoothness. In practice, the use of small parameters in  RBF-FD results in ill-conditioning, so that a stable algorithm is required for precise evaluation. Several approaches have been proposed to deal with the aforesaid ill-conditioning. An early study \cite{KansaHon2002} investigated six different approaches to deal with this problem. A first tool to work with multiquadric RBFs with small shape parameters was proposed by \cite{Wright2003,Forn2004} by removing the restriction that the shape parameter be real. This approach was called the Contour-Pad\'e algorithm and was limited to a small number of degrees of freedom only.

In 2007, \cite{Forn2007} proposed an RBF-QR approach to obtain a stable RBF algorithm for a special case of problems on the sphere. Later, \cite{Forn2011} applied similar ideas to obtain a stable algorithm for Gaussian RBFs in the Cartesian setting. Although RBF-QR provides excellent accuracy, its computational cost increases significantly with increasing shape parameter and in higher dimensions as well. Another approach for Gaussian kernels, called Gauss-QR, was proposed by \cite{Fass2012} by using the connection of the RBF-QR approach to Hilbert-Schmidt (or Mercer) series expansions of positive definite kernels. The RBF-GA algorithm was another development in this context, which does not require the use of truncated infinite series expansions or of any other form of numerical approximations \cite{Fornberg2013s}\footnote{This statement is only true if one accepts {\sc Matlab}'s implementation of the incomplete Gamma function as a method that does not approximate.}. Other modern methods include: weighted SVD \cite{DeMarchi20131} and Laurent series expansion approach \cite{Kindelan2016}. The most recent development in stable computation with `flat' RBFs is the RBF-RA algorithm using vector-valued rational approximation, which is basically an improvement of the Contour-Pad\'e algorithm \cite{Wright2017}.

An alternative approach is to use a shape parameter independent kernel by combining polyharmonic splines (PHSs) with augmented polynomials. \citet{Flyer2016} proposed that one can get a robust RBF-FD formulation by increasing the degree of polynomials higher than the order of PHS kernel. Following works \cite{Flyer20162,bayona2017role,fornberg2015primer} discussed, in detail, the role of polynomials in RBF-FD for interpolation and for solving PDEs, which was found to be much more than just to prevent singularities in 'non-favourable' node layouts. In the 'PHS+poly' approach, the order of convergence can end up higher than what the PHS part would suggest, since then the PHS coefficients tend to `naturally' vanish under node refinement, with the polynomial part `taking over' (being better able to fit smooth data). The PHS part still stabilizes the calculations, but does so without damaging the high accuracy levels provided by polynomials. Unlike the GA-based RBF-FD, the high condition number in the linear system is not a severe problem in "PHS+polynomial" case as it only reflects a theoretical limitation in the definition of condition number. The condition number in this may be sensitive to the matrix scaling without having any adverse effect on the computation. The stencil size has marginal effect on the accuracy and the order of convergence is practically independent of the order of PHS kernel. Therefore, the error convergence is entirely governed by the included polynomials, which results in a need for significantly more points in higher dimensions to achieve more accuracy. For example, with polynomial augmentation, a 12-point stencil is required to achieve second order accuracy. However, this could be obtained with only a 5-point stencil with standard RBF-FD \cite{Wright2006,Shankar2017} over a regular node-layout \cite{Shankar2017}.

Recently, a new approach to reduce the ill-conditioning problem in RBF approximations by using a hybrid Gaussian-cubic kernel was proposed \cite{Mishra2015AMC}. The basic idea behind such a hybridization is to obtain a kernel which utilizes the merits of two different kernels while compensating for the limitations of each and keeping the formulation as a standard RBF method. A follow-up work showed that such a hybrid kernel can improve the RBF-PS method for the numerical solution of PDEs for relatively large numbers of degrees of freedom using direct solvers and a reasonable compromise on the accuracy \cite{Mishra2017,MishraSEG2017}. In this paper, we present an RBF-FD formulation using the hybrid Gaussian-cubic RBF, which does not require the complicated machinery of the other stable algorithms (for Gaussian RBF-FD).

The rest of the paper is organized as follows. In section 2, we provide an introduction to RBFs in the context of solving PDEs. In section 3, we discuss the RBF-FD formulation inspired by \cite{Fass2015}. In section 4, we provide a numerical test, which shows the advantage of hybrid kernel-based RBF-FD over Gaussian RBF-FD, and how the hybrid kernel reduced the ill-conditioning in the resulting linear system. We also discuss the `stagnation error' and need for polynomial augmentation in the hybrid kernel and provide a heuristic comparison with 'PHS+polynomial' based RBF-FD. In section 5, we discuss the effect of kernel-parameters on the accuracy and stability of the RBF-FD formulation. Finally we provide some numerical examples for frequency-domain modeling of acoustic wave propagation in homogeneous and isotropic media.

\section{Radial kernels for PDEs}
The idea behind RBF interpolation was initially proposed for approximation of irregular surfaces on scattered datasets \cite{Hardy1971}, which was later established as an efficient interpolation scheme through a robust comparative study with many other approaches \cite{Franke1979}. Since then, many variants of RBFs have been introduced, theoretically developed, and applied to various problems in science and engineering including: scattered data interpolation, numerical solution of PDEs, computer graphics, surrogate modeling, and neural networks. In the context of numerical approximation of PDEs using RBFs, the choice of an appropriate kernel is crucial. Polyharmonic spline (PHS) kernels are piecewise smooth and provide relatively better conditioned linear systems than those obtained using many other standard RBFs. 
Table~\ref{tab:rbf} lists some typical radial basis functions. More information about various kernels including radial basis functions, in the context of numerical solution of PDEs can be found in \cite{Fassbook2007,Fass2015,Fornberg2015}.
\begin{table}
  \centering
  \begin{tabular}{|l|l|}
    \hline
    Kernel & Mathematical expression \\
    \hline
    Multiquadric (MQ)                       & $ (1+(\epsilon r)^2)^{1/2}$ \\
    Inverse multiquadric (IMQ)              & $ (1+(\epsilon r)^2)^{-1/2} $\\
    Gaussian (GA)                             & $ e^{-(\epsilon r)^2}$\\
    Polyharmonic Spline (PHS)                & $ \begin{cases} r^m \ln(r) \qquad m =2,4,6,... \\ r^m \qquad\qquad m= 1, 3, 5,... \end{cases}$ \\
    Wendland's (Compact Support)            & $(1-\epsilon r)^{4}_{+}(4\epsilon r+1)$\\
    Hybrid Gaussian-cubic                             & $\alpha e^{-(\epsilon r)^2}+\beta r^3$ \\
    \hline
\end{tabular}
 \caption{ Mathematical expressions of some radial kernels.}
 \label{tab:rbf}
\end{table}

The hybrid Gaussian-cubic RBF, proposed by \cite{Mishra2015AMC}, has been shown to provide stabilized schemes for scattered data interpolation problems. A follow-up work \cite{Mishra2017} has shown that  such a hybrid kernel is a reasonable choice for numerical solution of PDEs via the RBF-PS method and provides a stabilized RBF-PS method which works for a relatively large number of degrees of freedom. The hybrid Gaussian-cubic kernel is given as: 
\begin{equation}
\phi(r) = \alpha e^{-(\epsilon r)^2}+\beta r^3,
\label{eq:hkernel}
\end{equation}
\begin{figure}
\includegraphics[scale=0.4]{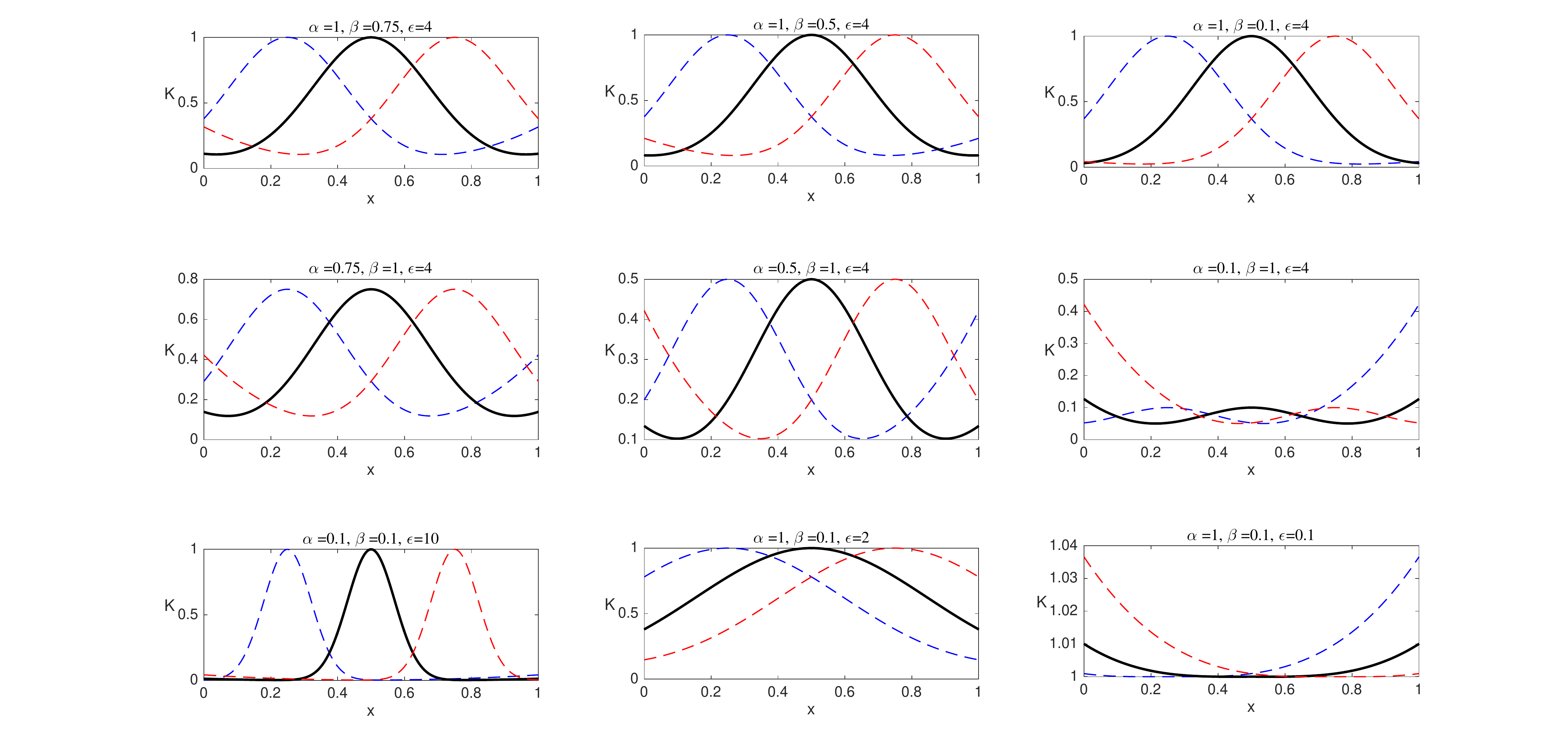}
\caption{ 1D visualization of the hybrid kernel, given by equation \ref{eq:hkernel}, for various parameter settings. (K = $\phi(r)$)}
\label{fig:kernelplot}
\end{figure}

\noindent where $\epsilon$ is the shape parameter for the Gaussian kernel and $\alpha$, $\beta$ are the weights controlling the contribution of the Gaussian and cubic kernels, respectively. The hybrid Gaussian-cubic kernel is visualized with various combinations of parameters in Figure~(\ref{fig:kernelplot}). Since scaling an RBF by a constant does not affect the algorithm, the above kernel can be normalized by introducing a factor $\gamma = \beta/\alpha$ as given by
\begin{equation}
\phi(r) = e^{-(\epsilon r)^2}+\gamma r^3.
\end{equation}

\noindent The new hybrid kernel has only two parameters, which shall control the accuracy and the stability of the RBF-FD algorithm presented here.

\section{RBF-FD Discretization}
RBF-FD is a local meshless method in which the discretization is obtained by computing a number of local differentiation matrices and assembling them into a single large, but sparse, system matrix. In this section, we take a general elliptic PDE and explain a step-by-step procedure to solve it by an RBF-FD scheme.

\noindent Let us assume a boundary value problem (Helmholtz type) in a rather general computational domain $\Omega$, as given by
\begin{equation}
\label{eq:general}
Lu(\bm{x}) = f(\bm{x}),
\end{equation}
where $L$ is a general linear differential operator ($\nabla^2+k^2$, for the Helmholtz equation), $u$ is a field, $f(\bm{x})$ is a source term and $\bm{x}$ is the spatial variable. A general boundary condition for equation (\ref{eq:general}) can be written as,
\begin{equation}
\label{eq:generalBC}
Bu(\bm{x}) = g(\bm{x}).
\end{equation}
Here $B$ is a general linear differential operator at the boundary. For Dirichlet boundary conditions $B$ is the identity operator and for Neumann boundary conditions $B=\frac{\partial}{\partial\bm{n}}$, where $\bm{n}$ is the unit outward normal at the corresponding boundary.

The first step in any meshless method is to create nodes inside the given computational domain. Figure~(\ref{fig:nodetype}a-\ref{fig:nodetype}d) shows some typical node arrangements, created using mathematical sequences, in a 2D domain. However, in principle, a meshless algorithm will work on any random distribution or manually defined node arrangements.
\begin{figure}
\begin{myfont}
\centering
\includegraphics[scale=0.45]{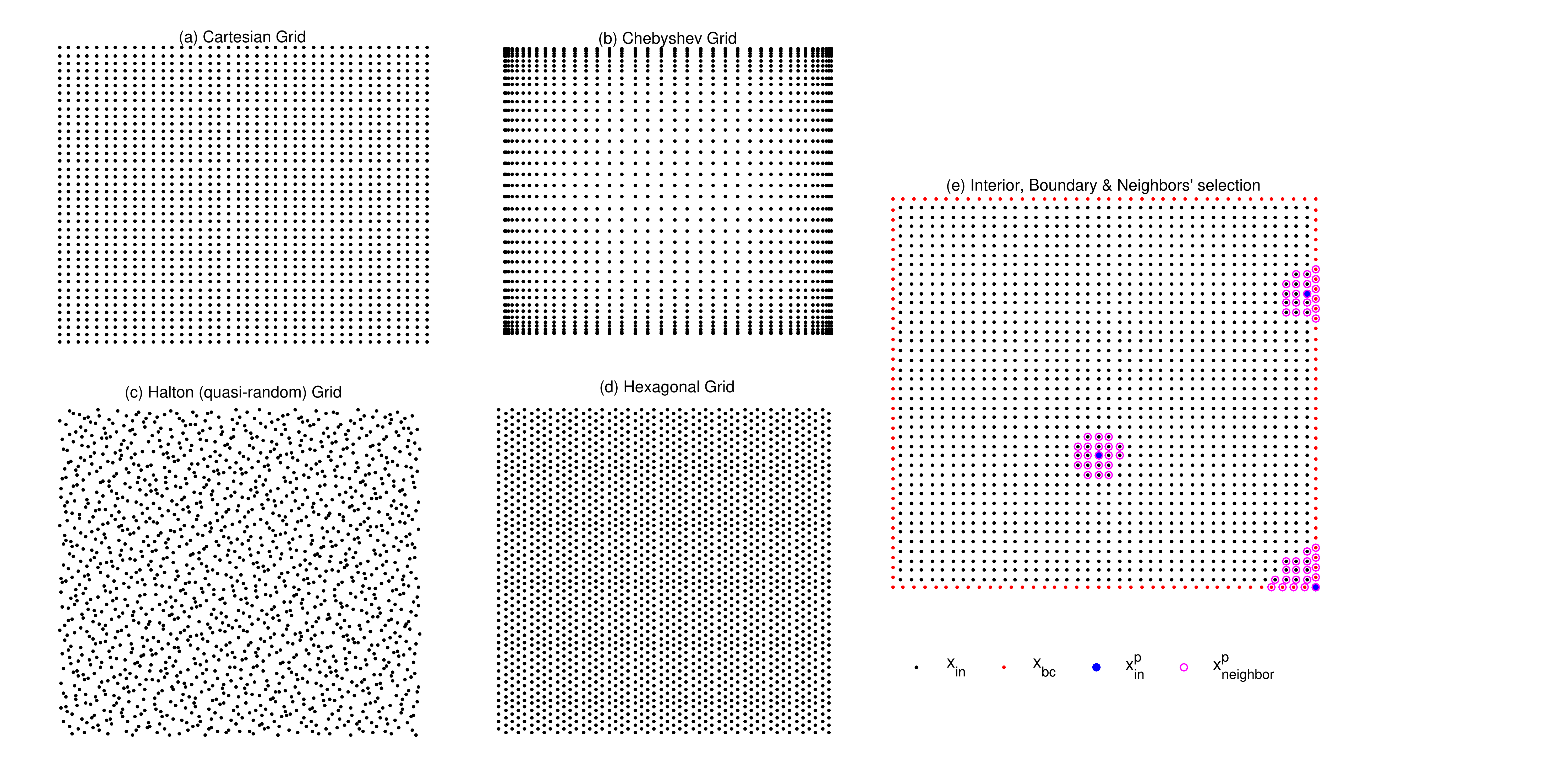}
\caption{ (a)-(d) Example of some typical node distributions in a 2D domain. (e) Visualization of the neighbors' selection for an interior, a boundary and a corner node, where 20 closest points ($x^{p}_{k}$) have been considered as neighbors for a $p^{th}$ node $x^{p}_{\text{in}}$. } \end{myfont}
\label{fig:nodetype}
\end{figure}

Since RBF-FD is a local meshless approach, the next step is to define a strategy for selection of a finite number of neighbor points for each node. In order to compute a differentiation matrix  at each point, only those neighbor points will be considered.  Selection of the neighbor points can be done by many approaches. One of them is to define a finite radius from a node, and consider the points inside the circle as neighbors of that node. This approach has been widely used for local meshless collocation methods and more recently for RBF-FD methods too. However, in order to maintain its accuracy near the boundary this strategy is often coupled with the use of `ghost nodes' outside the computational domain for computing the differentiation operator at boundary points \cite{Flyer2016}. In this paper, we choose a simple alternative for the selection of neighbors. Instead of defining a spatial support domain, for each node, we define a certain number of closest points as neighbors. Such an approach  for selection of neighbor nodes has been visualized in Figure (\ref{fig:nodetype}e).

In order to explain the discretization of a PDE (elliptic) via the kernel-based finite difference method, following \cite{Fass2015}, we start with the standard representation of interpolants.  An approximation of a field\footnote{By `field' we mean the unknown function in the PDEs, which needs to be approximated numerically. For example, in this work, the term `field' means the frequency-domain representation of the wavefield.} $u$, on the nodes $\bm{x}$ can be written using a specified radial kernel $\phi$ as
\begin{equation}
\centering
\hat{u}(\bm{x}) = \sum_{j=1}^{N} c_j \phi(\|\bm{x} - \bm{x}_j\|) =\bm{\Phi}(\bm{x})^T \mathbf{c}, \qquad \bm{x}\in \Omega,
\label{eq:inter}
\end{equation}
where $\bm{\Phi}(\bm{x})^T = \left(\phi(\|\bm{x} - \bm{x}_1\|,\phi(\|\bm{x} - \bm{x}_2\|),\ldots,\phi(\|\bm{x} - \bm{x}_N\|) \right)$, and $N$ is the number of nodes (in the domain $\Omega$) used for the approximation. We can also write equation (\ref{eq:inter}) in cardinal form as
\begin{equation}
\hat{u}(\bm{x}) = \bm{\Phi}(\bm{x})^T \mathbf{K}^{-1} \bm{u},
\end{equation}
where $\mathbf{K}$ is the global RBF interpolation matrix. The basic assumption here is that, since the kernel-based interpolation provides a good approximation $\hat{u}$ of a field $u$, any operator applied on $\hat{u}$ will be a reasonable approximation of the same operator applied on the true field \cite{Fassbook2007,Fass2015}. Also, it has been found that the hybrid kernel provides `good' approximation of the function \cite{Mishra2015AMC} as well as its derivative \cite{Mishra2017}. A linear operator $L$, applied on the field can therefore be approximated as
\begin{equation}
L \hat{u}(\bm{x}) = L\bm{\Phi}(\bm{x})^T \mathbf{K}^{-1} \bm{u}.
\end{equation}
If we perform such an approximation at the points $\bm{x}_1,\ldots,\bm{x}_N$ and use the notation
\begin{equation}
\mathbf{K}_L =
\begin{bmatrix}
L \bm{\Phi}(\bm{x}_1)^T \\
\vdots \\
L \bm{\Phi}(\bm{x}_N)^T
\end{bmatrix},
\end{equation}
then we can write the discretization $\mathbf{L}$ (an $N\times N$ matrix) of the linear operator $L$ as
\begin{equation}
\mathbf{L} = \mathbf{K}_{L} \mathbf{K}^{-1}.
\label{eq:globalD}
\end{equation}

At this stage, since the interpolation matrix is global,\emph{ i.e.}, all the points have been considered while computing the interpolation matrix, the discrete operator in equation (\ref{eq:globalD}) has a global nature. The idea for RBF-FD discretization is to consider only a `few' neighbor points to compute the discrete operator at a location (kernel centers). In order to have a more specific discussion using the neighbor nodes, let us assume that we have to compute the RBF-FD derivative at locations $\mathcal{X} = \{\bm{x}_1,\ldots,\bm{x}_N \}$. For the $i^{\text{th}}$ node $\bm{x}_i$, consider the number of neighbor nodes equal to $n_{\bm{x}_i}$. We will collect the stencils at locations $Z_i = \{\bm{z}_1,\ldots,\bm{z}_N\}$. Thus, the local differential operator at the stencil center $\bm{x}_i$ can be written as
\begin{equation}
\mathbf{L}_i = \mathbf{K}^{\bm{x}_i}_{L} \mathbf{K}^{-1}_{\bm{z}_i}.
\end{equation}

In an RBF-FD discretization, we compute various local differentiation matrices $(\mathbf{L}_i)$ and place them at specific locations in a global differentiation matrix $\mathbf{L}^{FD}$\footnote{This is a global (sparse) differentiation matrix in which all the local differentiation matrices are collected, therefore, it is different than the one in equation (\ref{eq:globalD})}, as given by
\begin{equation}
\mathbf{L}^{FD} =
\begin{bmatrix}
\mathbf{K}^{\bm{x}_1}_{L} \mathbf{K}^{-1}_{\bm{z}_1} \mathbf{P}_1 \\
\vdots \\
\mathbf{K}^{\bm{x}_N}_{L} \mathbf{K}^{-1}_{\bm{z}_N} \mathbf{P}_N \\
\end{bmatrix},
\end{equation}
where $\mathbf{P}_i \in \{0,1\}^{n_{\bm{x}_i}\times N}$ is an incidence matrix which has been defined to place the nodes $z_i$ (at which $\mathbf{L}_i $ has been computed) to the correct position in the sparse row of $\mathbf{L}^{FD}$. The elements of $\mathbf{P}_i$ are given by
\begin{equation}
[\mathbf{P}_i]_{k,l} = \begin{cases} 1 \qquad \text{if }k=l, \\ 0 \qquad \text{else.} \end{cases}
\end{equation}
Thus, the discrete representation of the problem defined by equations (\ref{eq:general}) and (\ref{eq:generalBC}), can be written as
\begin{equation}
\begin{bmatrix}
\mathbf{L}^{FD} \\
\mathbf{B}^{FD}
\end{bmatrix}[\mathbf{c}] =
\begin{bmatrix}
f \\
g
\end{bmatrix},
\label{eq:discrete}
\end{equation}
where $\mathbf{B}^{FD}$ is the discretized operator which needs to be applied at the boundary points, and calculated in a similar manner as $\mathbf{L}^{FD}$. Equation (\ref{eq:discrete}) has been written assuming only one kind of boundary condition, however, different boundary conditions can be incorporated by creating corresponding rows in it.

\section{Numerical Test 1}
We compare the RBF-FD formulations with four different types of RBF settings (1) Gaussian (GA) (2) Hybrid RBF, (3) PHS with high order polynomials and (4) GA with high order polynomials. The parameters for the hybrid kernel are set as $\epsilon=1$, and $\gamma = 10^{-6}$. In the 'PHS+polynomial' approach the degree of the augmented polynomial puts the constraint on the size of stencil, that is, to support a polynomial augmentation of degree $p$, the stencil size needs to be more than twice the number of the polynomial basis functions. For example, in 2D, if we want to augment a polynomial of degree 4, the stencil size should be at least 30 for stable formulation. In other words, the maximum allowed polynomial degree for a stencil size $30$, in 2D, would be 4. In this test, we fix the stencil size equal to 30 and compare the above mentioned four settings to formulate RBF-FD for computing the derivative of a simple function $f(x,y) = \sin(x)+\cos(y)$.
Figure (\ref{F:complot}) shows the error convergence (fill-distance (h) vs $l_{\infty}$- error) of different RBF-FD formulations, for computing the first-derivative of the above function. For a quasi-uniform node-layout, a usual estimate is  $h=\frac{(Vol(\Omega))^{1/d}}{N^{1/d}}$, where $N$ is total number of nodes in the domain $\Omega$. The stencil size is fixed to 30 and fourth order polynomials have been augmented. The slopes of the convergence plots have been estimated through a log-log fit on the whole data in each case. With GA RBF, for a fixed shape parameter, the convergence is disrupted. On the other hand, the hybrid kernel with fixed shape parameter shows a good convergence of roughly a similar order as the 'PHS+Poly' approach. However, if we add the same polynomial to the hybrid RBF, we could end up achieving even better accuracy with the hybrid RBF. The reason for this is that the major part of the hybrid RBF is GA, which is infinitely smooth and hence corresponds to a better accuracy. 

\begin{figure}
\centering
\includegraphics[scale=0.6]{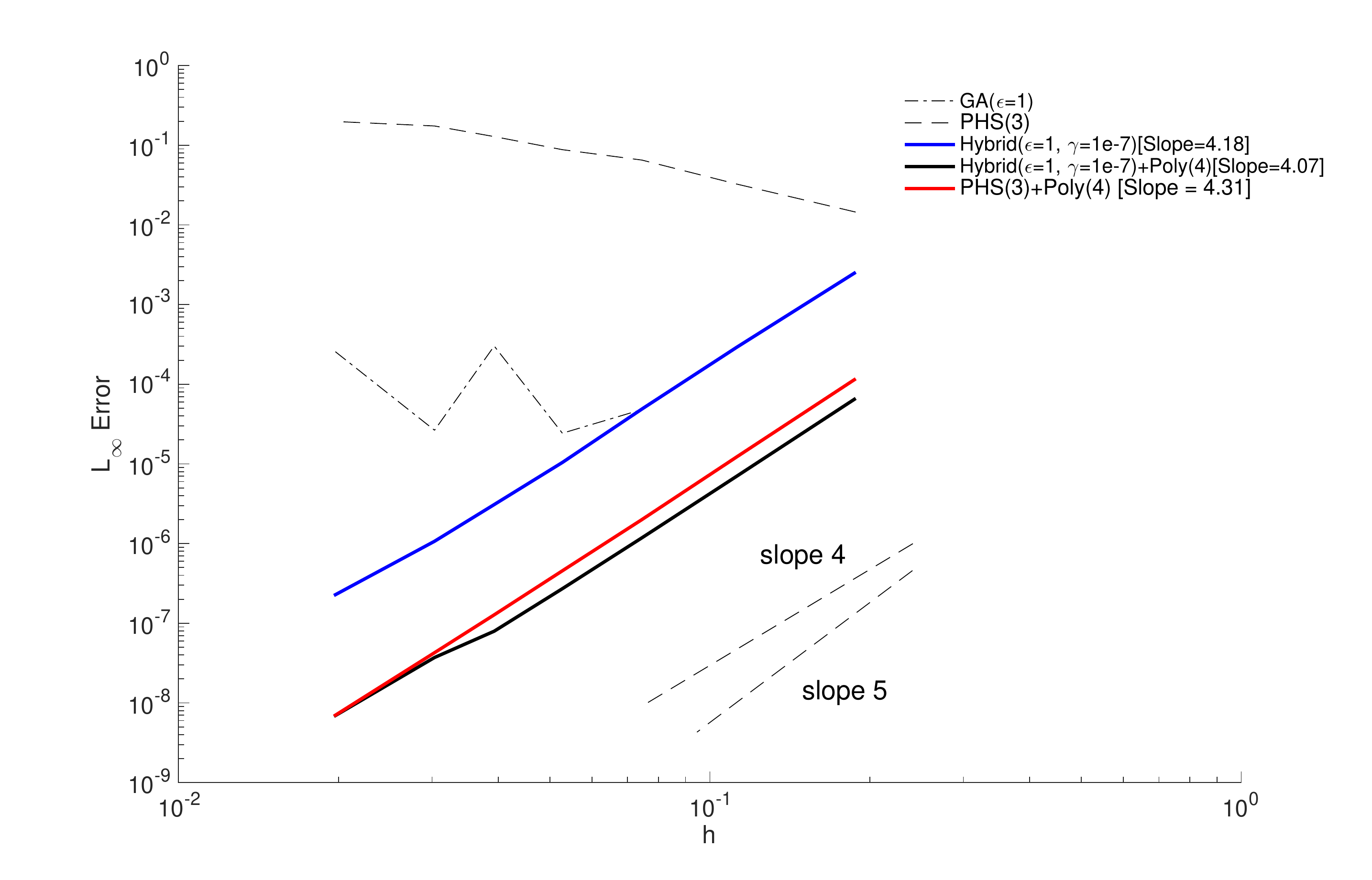}
\caption{\small Comparison of convergence plots for first derivative approximation with RBF-FD with different kernels: Gaussian (GA), third order polyharmonic spline (PHS(3)), and the hybrid kernel. Poly(4) stands for a fourth order polynomial augmentation.}
\label{F:complot}
\end{figure}

\section{Numerical Test 2}
For this test, we consider a typical 2D-PDE coupled with Dirichlet and Neumann boundary conditions as given by
\begin{equation}
\label{eq:helmholtz}
\nabla^2 u(x,z) - k^2 u(x,z) = 2\cos(x^2+z)-(4x^2+1+k^2)\sin(x^2+z), \qquad (x,z) \in \Omega, 
\end{equation}
\begin{equation}
u(x,z)  = \sin(x^2+z), \qquad (x,z) \in \Gamma -\Gamma_4,
\end{equation}
\begin{equation}
\frac{\partial u(x,z)}{\partial \bm{n}}  = \cos(x^2+z), \qquad (x,z) \in \Gamma_4,
\end{equation}
where $\Omega = [-1,1]\times[-1,1]$ and $\Gamma_4$ is the boundary having $z=1$. The analytical solution of this problem is
\begin{equation}
u(x,z) = \sin(x^2+z).
\end{equation}

First of all, we test the sparsity pattern and eigenvalue stability of the system matrix for this problem using the parameters $k=9$, $N=400$, $n=10$, $\epsilon=1$, and $\gamma=0.001$. According to \cite{flyer2014radial}, a natural intrinsic  irregularity in the stencil causes eigenvalues of the differentiation matrix to scatter into the right-half of the spectrum. This makes the RBF-FD algorithm unstable, for solving dissipation-free PDEs. Moreover, as the stencil size increases, so does the scatter of the eigenvalues in the right-half of the complex plane, which may be too large to handle even PDEs with dissipation. The eigenvalue-based instability remains a concern even in 'PHS+poly' based RBF-FD as the polynomial degree increases and so does the stencil size --- as discussed in \cite{barnett2015robust}, which needs special stabilization approaches. Figure (\ref{fig:eigen_distri}) shows the sparsity and eigenvalue spectrum of the differentiation matrix with GA and hybrid RBF for (1) evenly-spaced Cartesian and (2) random node distributions. It was observed that, irrespective of the irregularity in the stencils, the eigenvalues of the differentiation matrix remain stable when using the hybrid kernel. The effect of increasing stencil size on the eigenvalue stability has been shown in the figure (\ref{fig:eigen_num}). The eigenvalues in hybrid kernel-based RBF-FD remain stable even if the stencil size increases unlike the eigenvalues of GA-based RBF-FD.

\begin{figure}[t]
\flushleft
\includegraphics[scale=0.85]{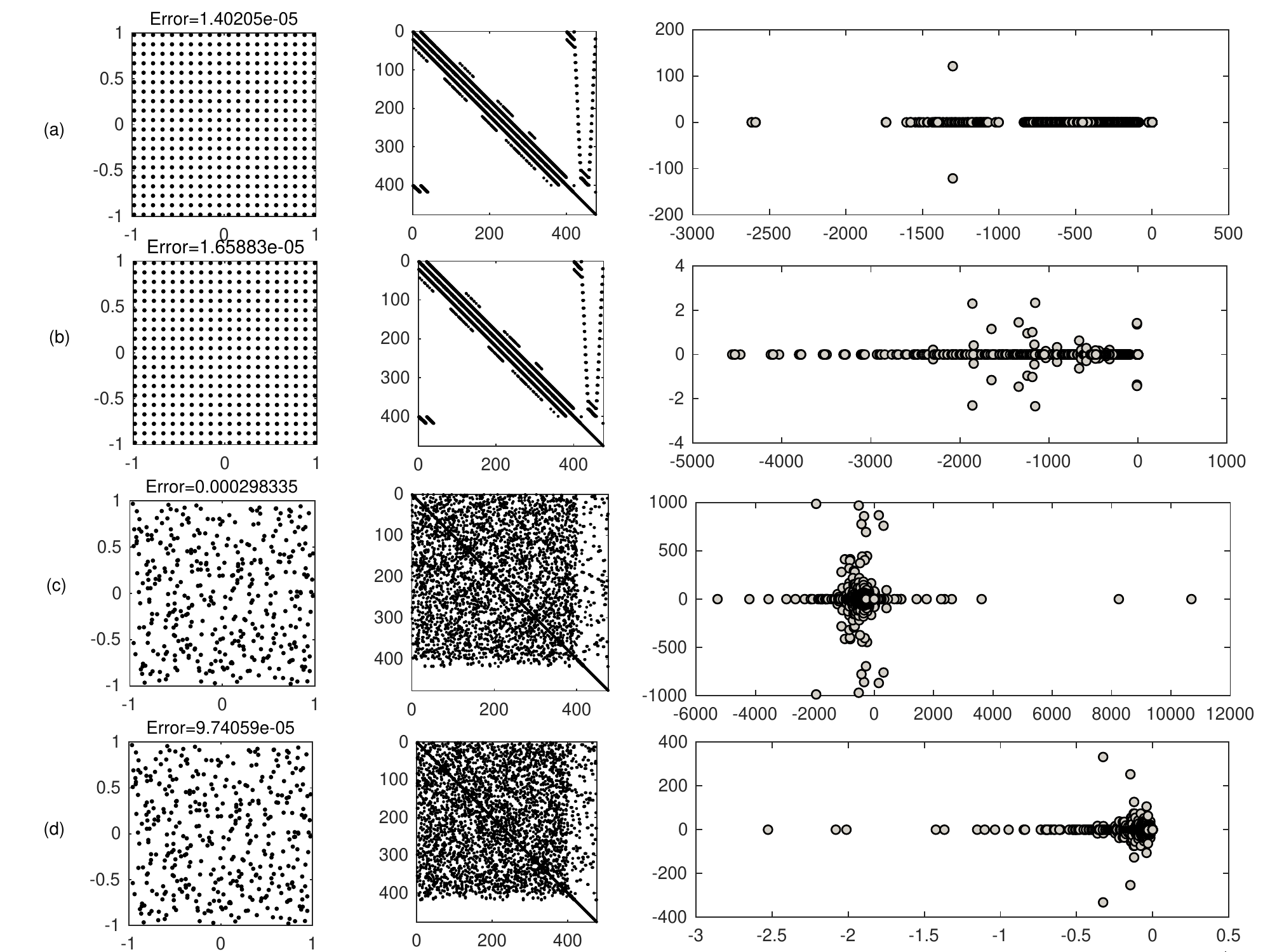}
\caption{ \small Sparsity and eigenvalue spectra of the differentiation matrix with GA (a, c) and hybrid RBF (b, d) for (1) evenly-spaced Cartesian and (2) random node distributions. }
\label{fig:eigen_distri}
\end{figure}

\begin{figure} [t]
\flushleft
\includegraphics[scale=0.85]{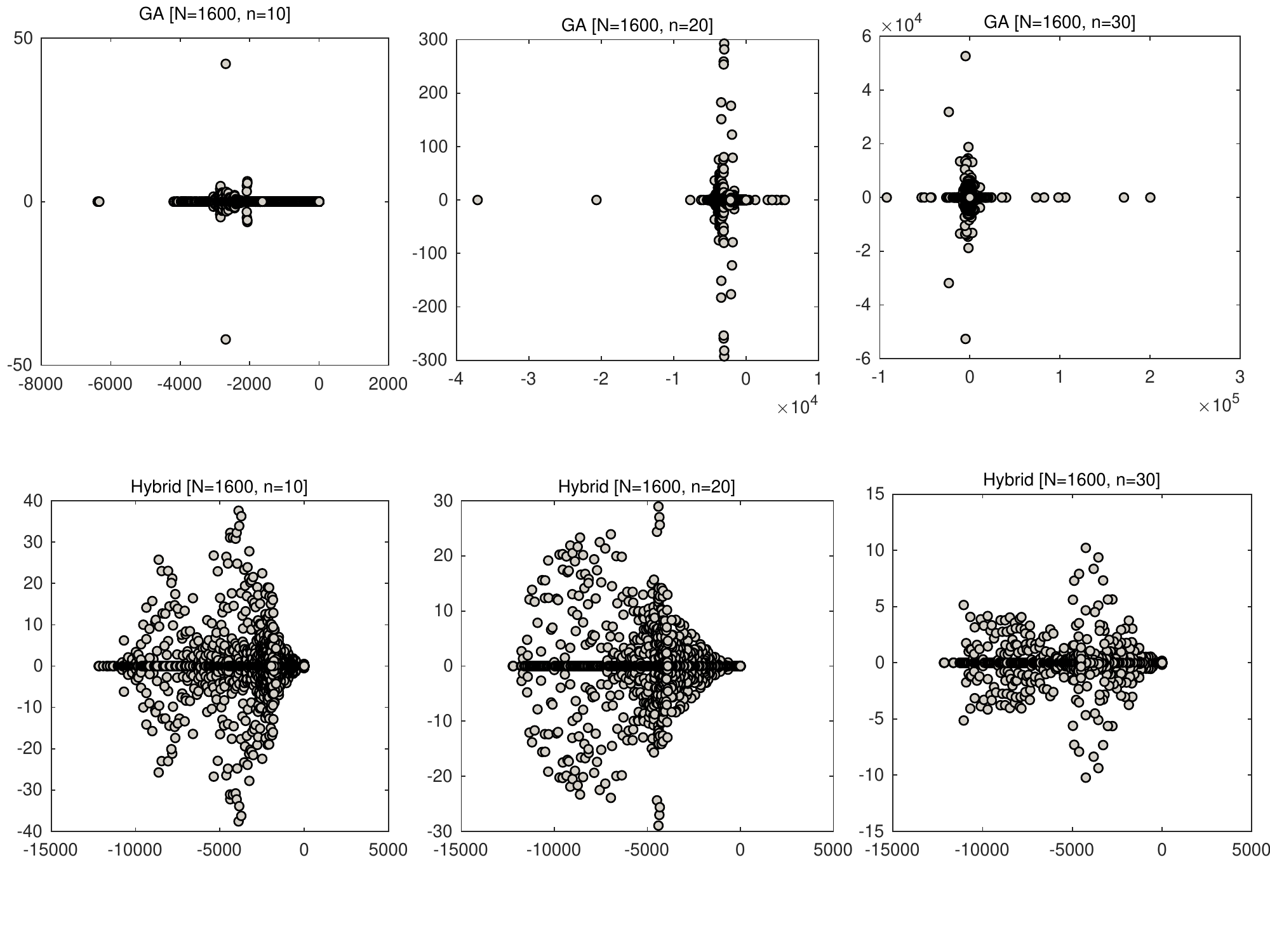}
\caption{ \small Eigenvalue spectra of the differentiation matrix with GA and hybrid RBF for different stencil sizes. $N$ is total number of nodes in the domain and $n$ is the stencil size. }
\label{fig:eigen_num}
\end{figure}

We further investigate the effect of kernel parameters on the efficacy of the algorithm. Figure~(\ref{fig:converg}) shows a contour visualization of error in the solution with various combinations of the shape parameter ($\epsilon$) and the weight ($\gamma$).  Whereas, the variation in the condition numbers of the interpolation matrix, for various combinations of $\epsilon$ and $\gamma$, have been visualized in Figure~(\ref{fig:Cwith}). In order to understand the advantage of using the proposed hybrid kernel in an RBF-FD algorithm, Figures~(\ref{fig:converg}) and (\ref{fig:Cwith}) should be interpreted together.  

\begin{figure}
\flushleft
\includegraphics[scale=0.45]{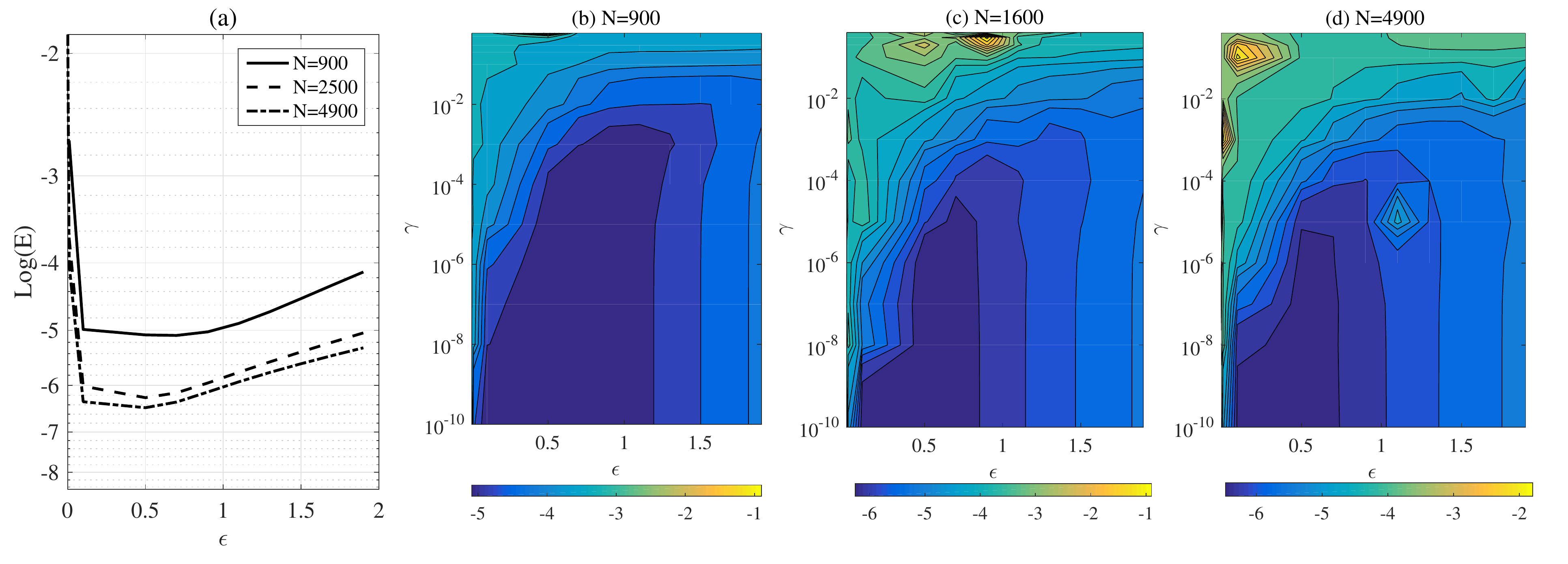}
\caption{\small (a) Error variation for a range of $\epsilon$ values for different numbers of degrees of freedom. This plot shows that pure Gaussian RBF ($\gamma =0)$ leads to a loss in accuracy for small shape parameters. Interpretation of this plot together with Figure~(\ref{fig:Cwith}a) suggests that using only Gaussian kernels in RBF-FD leads to an unstable formulation. (b)-(d) Error ($\log_{10}(\text{E})$) contours for various combinations of $\epsilon$ and $\gamma$ for different numbers of degrees of freedom. Number of neighbors is equal to 10. }
\label{fig:converg}
\end{figure}

\begin{figure}
\flushleft
\includegraphics[scale=0.45]{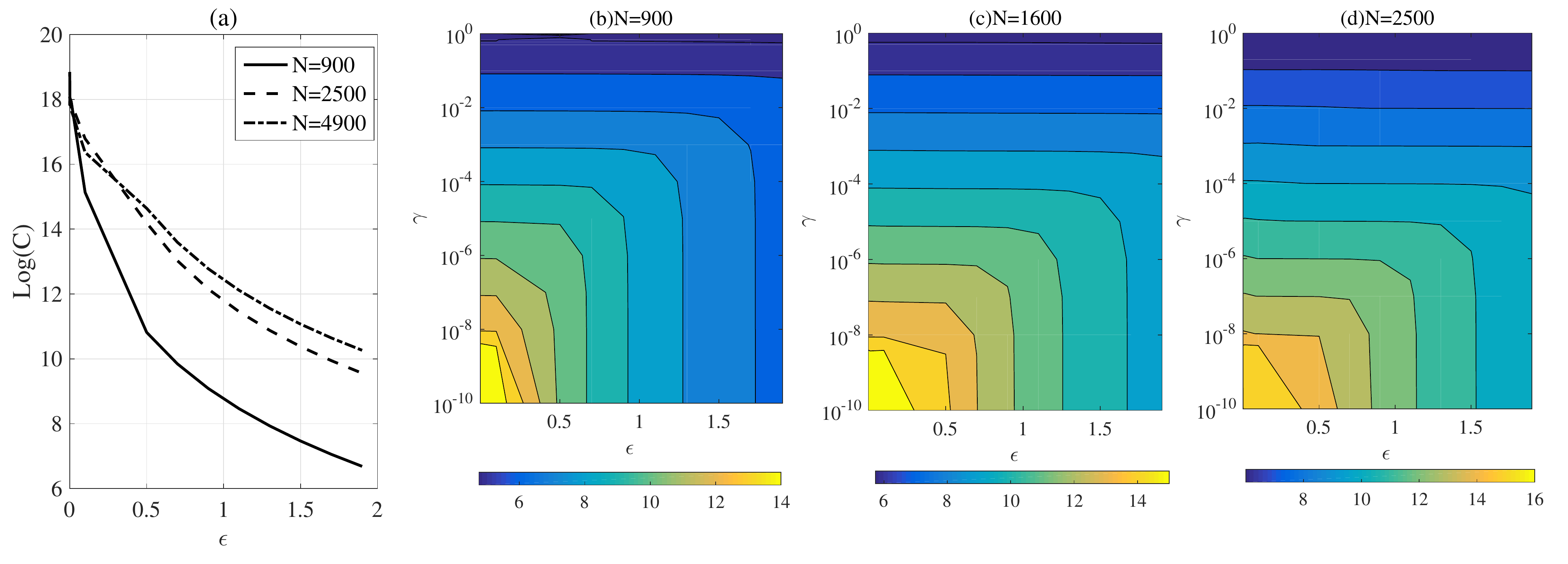}
\caption{ \small (a) Condition number variation for a range of $\epsilon$ values for different numbers of degrees of freedom. This plot shows that pure Gaussian RBF ($\gamma =0$) leads to ill-conditioned interpolants for small shape parameters as well as for higher numbers of degrees of freedom. Contours of condition numbers ($\log_{10}(\text{C})$) of system matrices for various combinations of $\epsilon$ and $\gamma$ at different numbers of degrees of freedom. Number of neighbors is equal to 10.}
\label{fig:Cwith}
\end{figure}

In spite of being a local approximation, application of pure Gaussian RBF in an RBF-FD algorithm may lead to a linear system with ill-conditioned interpolant as shown in Figure~(\ref{fig:Cwith}a), which corresponds to a significant loss in the accuracy as shown in Figure~(\ref{fig:converg}a). Such an instability is either due to the use of a small value of the shape parameter or due to an increase in the problem size, i.e., number of degrees of freedom. Such an ill-conditioning problem can be taken care of by using the hybrid Gaussian-cubic RBF as shown in Figures~(\ref{fig:converg}) and (\ref{fig:Cwith}). These observations are of particular interest in the context of RBF-FD because the RBF-FD algorithm with Gaussian kernels being severely ill-conditioned  requires computationally expensive algorithms like RBF-QR or RBF-GA to solve it. However, application of the hybrid kernel makes the system well-posed, such that is can be solved using a direct method, which in turn makes the algorithm significantly faster. We will discuss the computational cost of the present approach later in this paper. It was  observed that the error plots suggest smaller values of $\gamma$ for better accuracy, whereas the condition plots suggest the larger values of $\gamma$ for better stability. 
This can be interpreted as an instance of the well-known ``trade-off phenomenon" or ``uncertainty principle" in kernel-based approximation \cite{schaback2007practical,Fass2015}. The general idea behind this principle is that one cannot simultaneously achieve good conditioning and high accuracy using the ``standard basis".  The tension between numerical stability and accuracy may be viewed from different perspectives. To obtain a more stable formulation our options essentially are (1) to find a ``better basis"\citep{beatson2001fast,Forn2011,Fass2012}, or (2) to find ways to work with an easily computable "standard"  basis. Our hybrid kernel approach corresponds to a change in the function space in which we inherit good stability from the cubic kernel as well as high accuracy from the Gaussian kernel, and thus falls under option (2). 

 Inferring from the previous tests, we select the parameters of the hybrid kernel as $\epsilon =0.9$ and $\gamma =0.001$ and test the convergence of the presented RBF-FD approach for this case as shown in Figure~(\ref{fig:hybridvscubic}). 

\begin{figure}[t]
\centering
\includegraphics[scale=0.45]{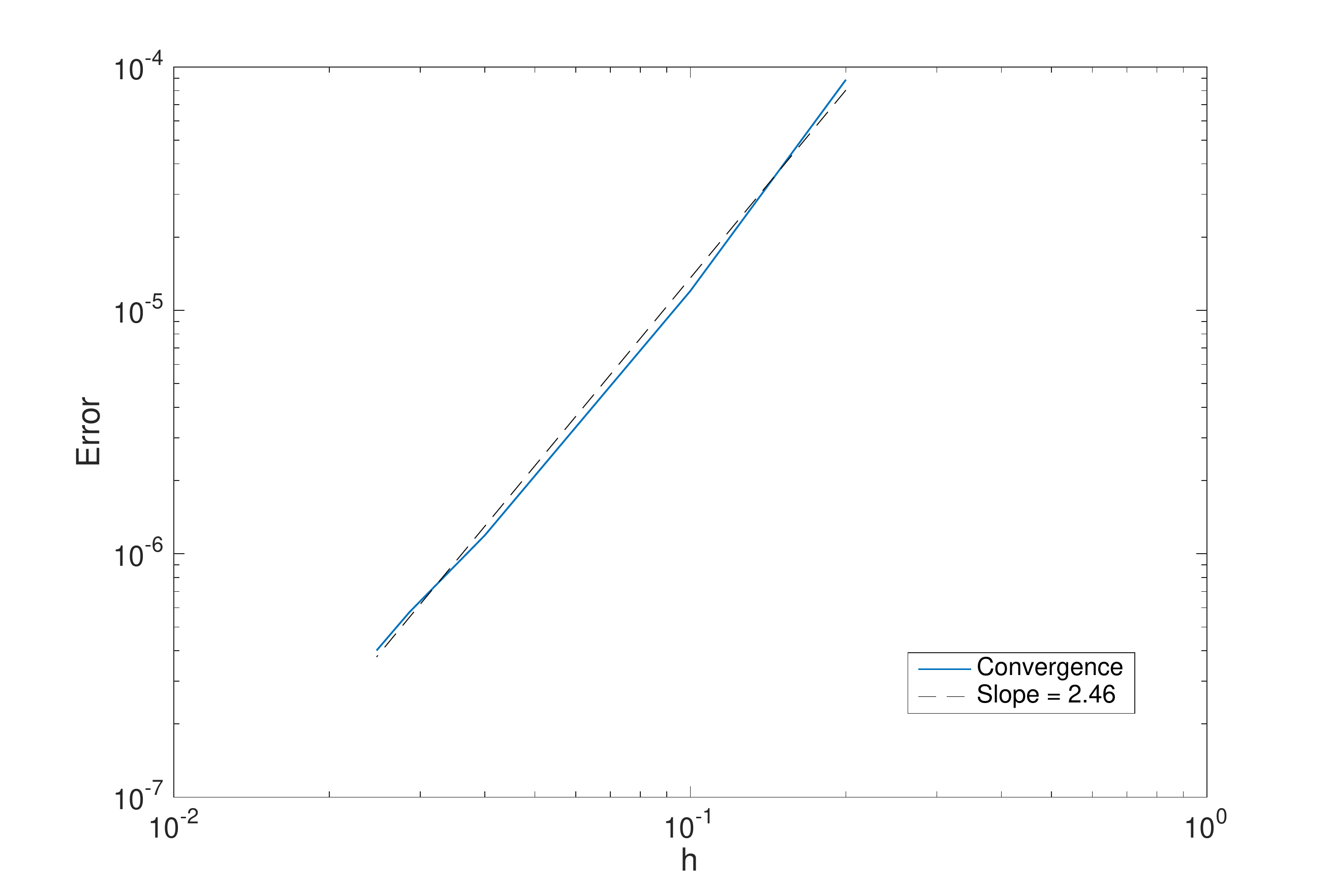}
\caption{Error convergence of RBF-FD with hybrid kernel corresponding to the numerical test in this section. Number of neighbors is fixed ($=10$).}
\label{fig:hybridvscubic}
\end{figure}

\begin{figure}
\centering
\includegraphics[scale=0.45]{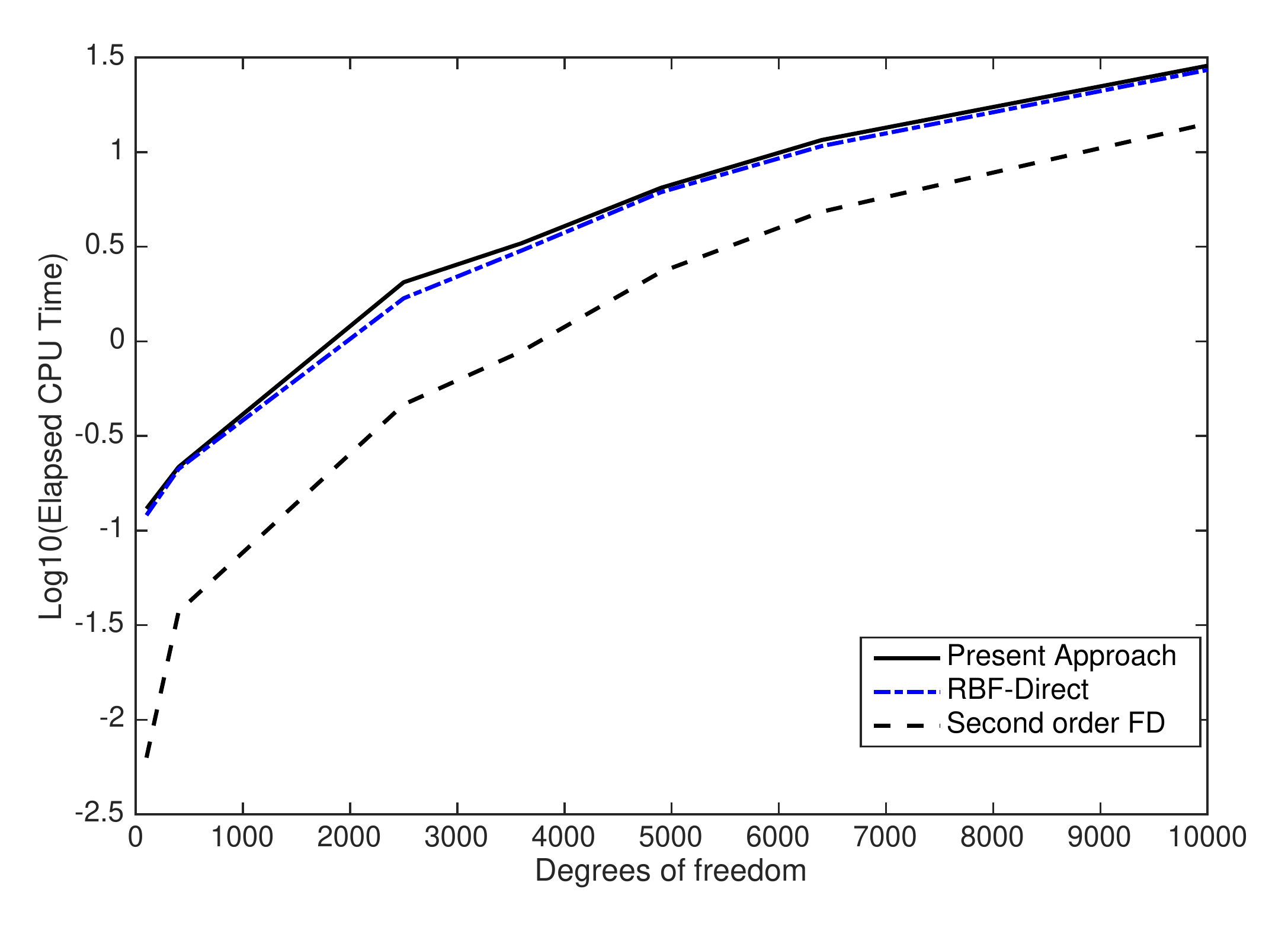}
\caption{Elapsed CPU time comparisons. RBF-FD with hybrid kernel corresponds to the numerical test used for the convergence plot with same parameters. Number of neighbors is fixed ($=10$).}
\label{fig:ctime}
\end{figure}

\subsection{\large Computational Time}
\noindent  Here, we provide the computational cost of the RBF-FD algorithm with hybrid kernel. We use the previous numerical test with the same parameters for this. The application of the hybrid kernel makes the system well-posed so that it can be solved using a direct method, which in turn makes the algorithm significantly faster. Figure~(\ref{fig:ctime}) shows the elapsed CPU time for solving the PDE in numerical test 2. Since using the hybrid kernel does not change the underlying structure of the direct version of RBF-FD, the computational cost of the present approach is practically the same as for the RBF-Direct method with the Gaussian RBF. This suggests that using the hybrid kernel in the RBF-FD method leads to a fast and stable algorithm. For completeness, we also include the computational cost of a second order finite difference method for the same test. However, it should be noted that stable methods are often more accurate than RBF-direct approaches like 'PHS+poly' or the present one. \footnote{The state-of-the-art is to provide a comparison of the computational cost of a new approach against the RBF-Direct method such that a heuristic comparison with other available approaches can be made without actually using them, as done by \cite{Wright2017}.} 

\section{Numerical Test 3}
Frequency-domain modeling of acoustic wave propagation involves numerical approximation of the Helmholtz equation, for which the finite difference method has been a convenient choice \cite{Dablain1986,Charl1996,Changsoo1998,Hustedt2004,Etgen2007,Amini2011,TaoSen2013,Moreira2014,Liu2014,TakekawaT145}.  In order to suppress the spurious reflections from the truncated computational boundary, absorbing boundary conditions are coupled with the governing Helmholtz equation \cite{Clayton1977,Cerjan1985,Berenger1994,Liu2010}.

Now that we have established the stability and convergence of the proposed RBF-FD approach, in this section, we perform simple numerical tests using RBF-FD in the context of frequency-domain solution of wave propagation in homogeneous and isotropic media as a preliminary application. For all the cases, we have kept the number of neighbor nodes equal to $10$, and the parameters of the hybrid kernel as $ (\epsilon=1$ and $\gamma =10\text{e-06})$.
 
In a Cartesian coordinate system, the 2D constant-density time-domain acoustic wave equation is given as
\begin{equation}
\nabla^2 p(\mathbf{x},t)-\frac{1}{c(\mathbf{x})^2}\frac{\partial^2 p(\mathbf{x},t)}{\partial t^2}=f(\mathbf{x},t),
\end{equation}
where $c$ is the primary wave velocity, $p(\mathbf{x},t)$ is the pressure wavefield, $f(\mathbf{x},t)$ is the source term, $\mathbf{x} = (x,z)$ are spatial coordinates, and $\nabla^2$ is the Laplacian operator. However, in the context of full waveform inversion and modeling the attenuation process, it is often required to solve the acoustic wave equation in frequency-domain, which is given by
\begin{equation}
\label{eq:seismic}
\nabla^2 \tilde{p}(\mathbf{x},\omega)+\frac{\omega^2}{c(\mathbf{x})^2}\tilde{p}(\mathbf{x},\omega) = \tilde{f}(\mathbf{x},\omega),
\end{equation}
where $\omega$ is angular frequency, and $\tilde{p}(\mathbf{x},\omega)$ and $\tilde{f}(\mathbf{x},\omega)$ are frequency-domain wavefield and source term, respectively, as given by
\begin{equation}
\tilde{p}(\mathbf{x},\omega) = \frac{1}{\sqrt{2\pi}} \int_{-\infty}^{+\infty} p(\mathbf{x},t) e^{-i\omega t} dt,
\end{equation}
\begin{equation}
\tilde{f}(\mathbf{x},\omega) = \frac{1}{\sqrt{2\pi}} \int_{-\infty}^{+\infty} f(\mathbf{x},t) e^{-i\omega t} dt.
\end{equation}

In order to suppress the spurious reflections from the truncated computational boundary, we include absorbing boundary conditions as given by
\begin{equation}
\label{eq:seismicbc}
\frac{\partial \tilde{p}(\mathbf{x},\omega)}{\partial \bm{n}} + i\frac{\omega}{c(\mathbf{x})}\tilde{p}(\mathbf{x},\omega)=0.
\end{equation}

\subsection{\large Frequency-domain visualization}
We start with the simple test of solving the above frequency-domain problem in the spatial domain $[0,1]\times[0,1]$ with a constant velocity $c= 1$. The spatial domain is discretized using $60\times60$ equally spaced Cartesian nodes. We incorporate the discrete Dirac-delta function as the energy source:
\begin{equation}
s(x,z, \omega) = \frac{1}{(h_xh_z)^2}\delta(sx)\delta(sz)
\end{equation}
for a point source located at $(s_x, s_z) = (0.2, 0.8)$, and $(h_x, h_z)$ are node spacing in the corresponding dimensions. Figure~(\ref{fig:exactRBFFDFD}) shows the approximate solution of this problem using RBF-FD and its comparison to the exact solution, which are in good agreement with the root mean square error equal to $9.2\text{e-03}$. This accuracy is less than the one we saw in the previous numerical test, which is due to the `imperfection' in the used absorbing boundary conditions \cite{Clayton1977}, which although it minimizes the spurious reflections from the computational boundary, it does not completely remove them. Figure (\ref{fig:NT2error}) shows the comparison of amplitudes (1) exact (2) computed by RBF-FD method with hybrid kernel.
\begin{figure}
\centering
\includegraphics[scale=0.45]{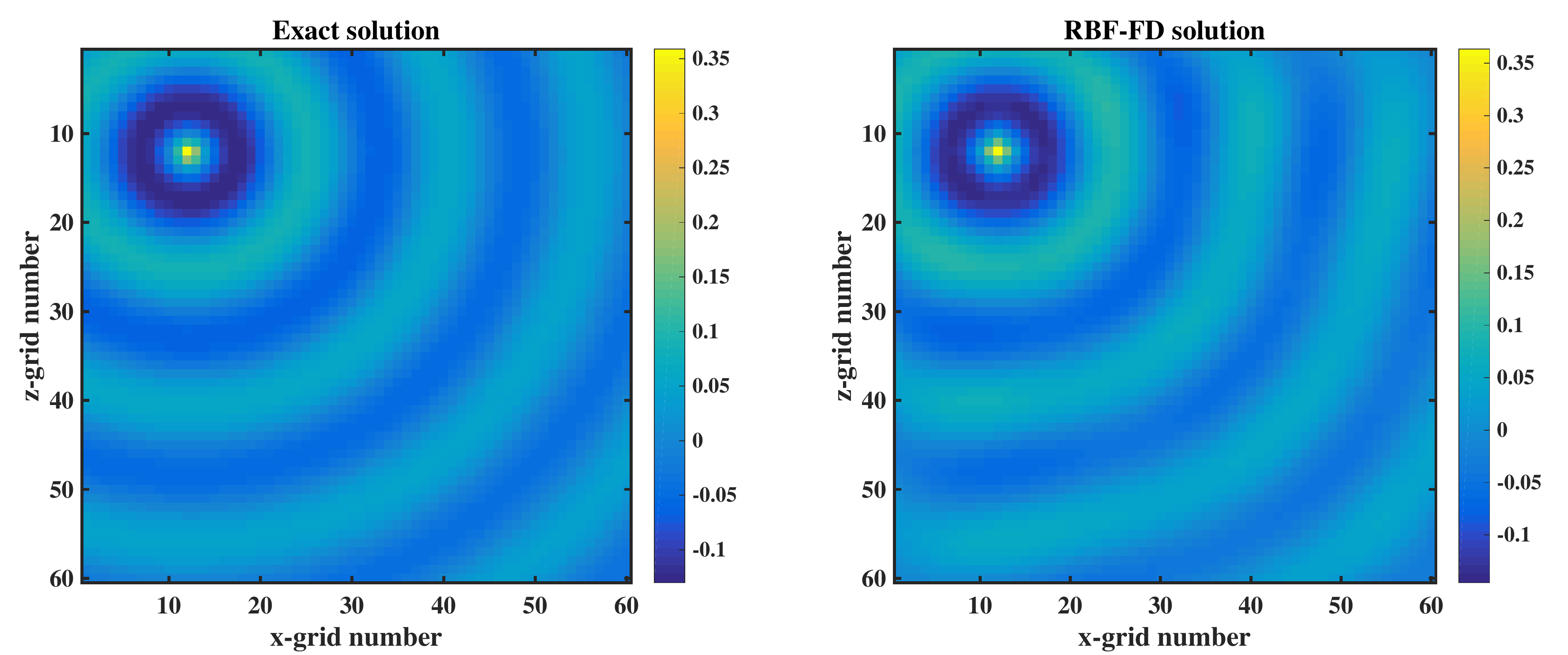}
\caption{ Approximate solution of the Helmholtz problem for a Dirac delta source, using RBF-FD compared with the analytical solution, \textit{i.e.}, Hankel function.}
\label{fig:exactRBFFDFD}
\end{figure}

\begin{figure}
\hspace{-1cm}
\includegraphics[scale=0.6]{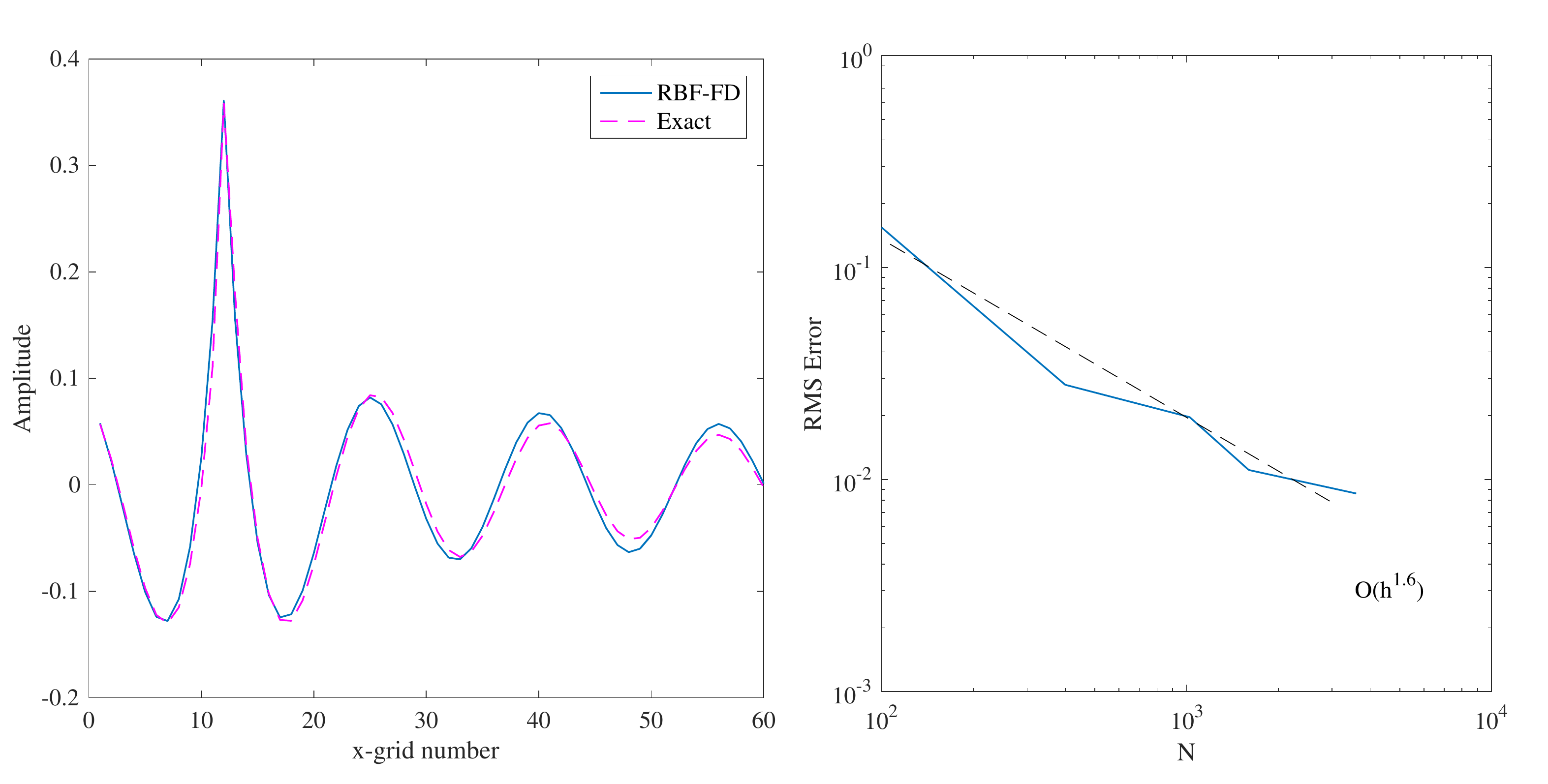}
\caption{(a) Comparing amplitude of wavefield passing through the source. (b) Error convergence with increasing degrees of freedom.}
\label{fig:NT2error}
\end{figure}

Now we solve the similar problem in a larger domain, \textit{i.e.}, $[0m, 400m]\times[0m, 400m]$. The aim of this test is to examine the solution at a practical range of frequencies, which is used for acoustic forward modeling. We take a typical Ricker source, which in the frequency domain is given by
\begin{equation}
\label{eq:RickerS}
  s(x,z,f) =  \frac{2f^2}{\pi^2 f^{3}_{c}} \exp \left( - \frac{f^2}{\pi f^2} \right) \delta(s_x)\delta(s_z),
\end{equation}
where $f = \omega / {2\pi}$ and $f_c = \bar{f}(3\sqrt{\pi}$) is related to the cuttoff frequency $\bar{f}$ \cite{Moreira2014}. We solve the problem for three different frequencies, $10$, $25$, and $50$ Hz, by keeping the total number of nodes in the domain as constant ($50\times50$). Figure~(\ref{fig:homo4freq}) shows the frequency-domain approximated wavefield for aforementioned frequencies, which suggests that the approximated solution does not disperse at relatively high frequencies.

\begin{figure}
\hspace{-1cm}
\includegraphics[scale=0.30]{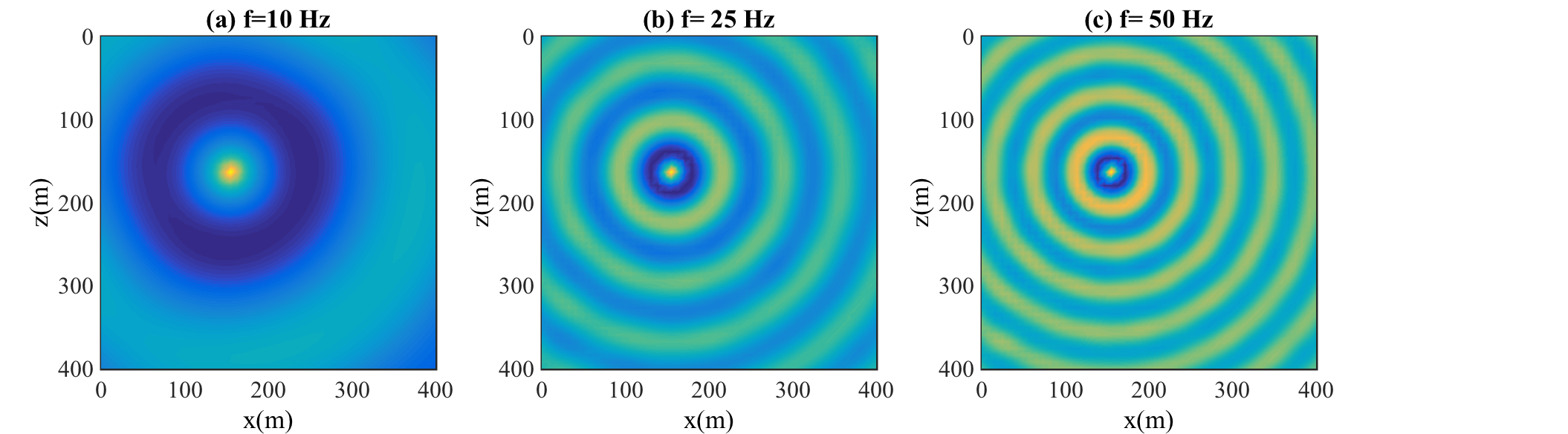}
\caption{Solution of the frequency-domain acoustic problem with Ricker source at frequencies 10, 25, and 50 Hz. The domain is discretized using $50\times50$ regularly spaced nodes for all three cases. The velocity in this medium has been taken as 2000m/s.}
\label{fig:homo4freq}
\end{figure}

In numerical test 1, it was shown that the eigenvalues of the `system matrix' remain stable for all kind of node arrangements. Figure~(\ref{fig:IRRgrid}) shows the solution of the frequency-domain wave equation with a 10Hz Ricker source on quasi-random nodes, generated through the Halton sequence, compared to the solution obtained on evenly spaced Cartesian grids. The two solutions are in reasonable agreement. This illustrates the meshless nature of the RBF-FD method.

\begin{figure}
\centering
\includegraphics[scale=0.3]{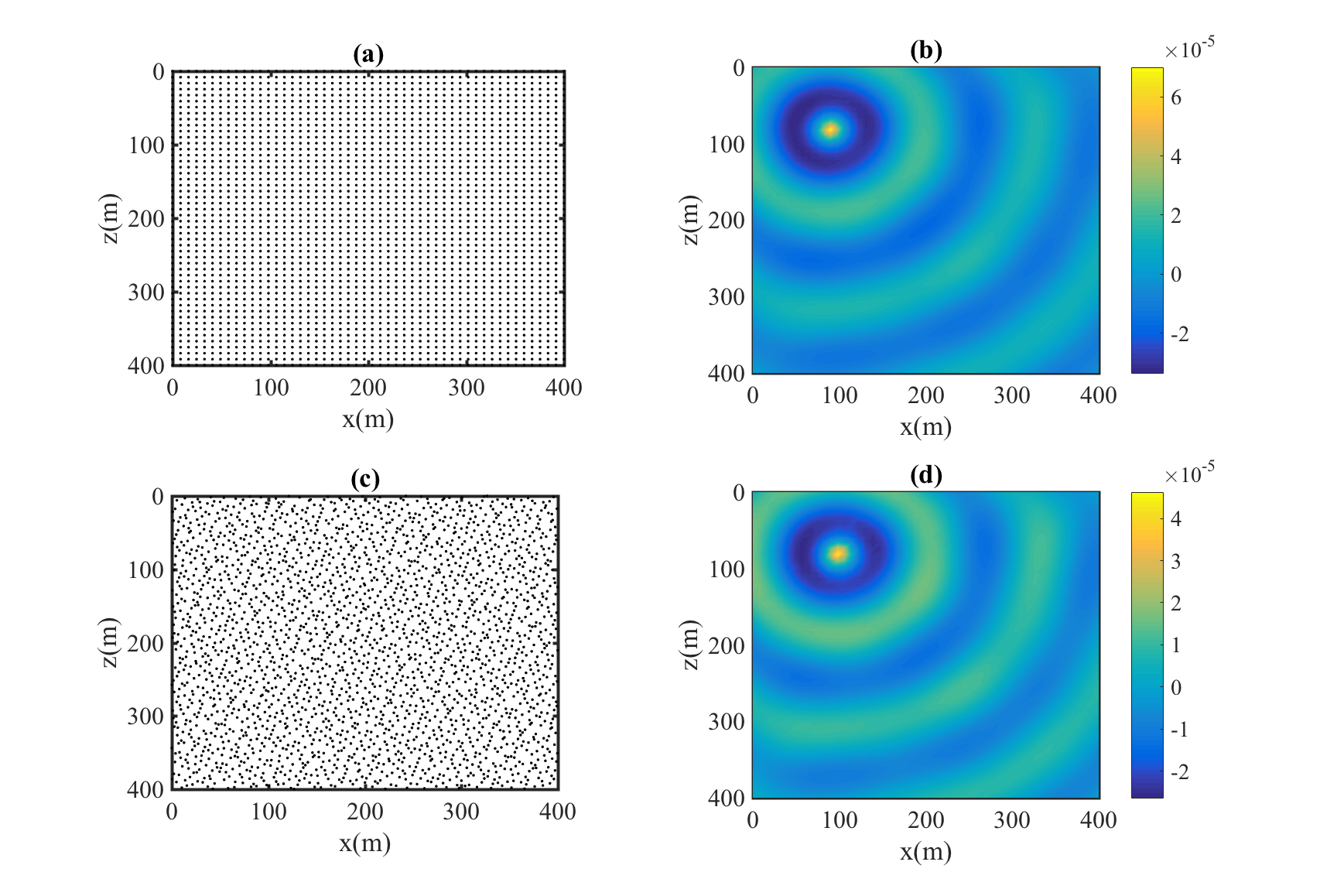}
\caption{Solution of the frequency-domain acoustic problem with Ricker source at 10Hz frequency. The domain is discretized using $50\times50$ (a) evenly spaced nodes and (c) quasi-random (Halton) nodes, and the corresponding solution of the wave equations is shown in (b) and (d), respectively. The velocity in this medium is taken as 2000m/s.}
\label{fig:IRRgrid}
\end{figure}

\subsection{Time-domain visualization}
In this section, we solve the acoustic forward problem given by equations (\ref{eq:seismic}) and (\ref{eq:seismicbc}) in the frequency domain for a set of frequencies and then transform the solution into the time domain by using the inverse Fourier transform. We compare the time-domain RBF-FD solution with that obtained by a 9-point mixed-grid finite difference method. We compute the FD solution by following the forth-order mixed grid formulation, and a different Ricker source (with both the FD and the RBF-FD) used in \cite{Amini2011}, which is given by
\begin{equation}
s(x,z,f) = \sqrt{\frac{4}{\pi(f_0)^2}} \left(\frac{f}{f_0}\right)^2 \exp{  \left(\frac{f}{f_0}\right)^2} \exp{\left(\frac{-2i\pi f}{f_0}\right)},
\end{equation}
where $f_0$ is the dominant frequency of the Ricker source. We consider a homogeneous medium having dimensions $[0m, 300m]\times[0m, 300m]$ and primary wave velocity as $2000m/s$. The domain is discretized using $900$ interior nodes (evenly spaced), having node interval equal to $10m$. The Ricker source is located at $(150m, 10m)$. An array of receivers is deployed at the same depth as the source and uniformly distributed between $x=10m$ to $x=290m$ as shown in Figure~(\ref{fig:seismo}a). The problem is solved in the frequency domain at a set of frequencies between 0 $Hz$ and 80 $Hz$, having 1 $Hz$ frequency intervals. The frequency-domain solution is then transformed into time-domain with a sample interval equal to $0.0125s$ for a maximum time of $1s$. The acoustic wavefield at $t= 0.1750s$ is shown in Figure~(\ref{fig:seismo}b), whereas Figure~(\ref{fig:seismo}c) shows the shot gather for this test. A comparison of the middle trace, passing through the source is displayed in Figure~(\ref{fig:seismo}d) and compared with that computed by using the 9-point mixed grid FD method with $100m$ thick perfectly matched layers.

\section{Conclusion}
In this paper, we have presented a stabilized RBF-FD formulation by using the hybrid Gaussian-cubic kernel and its application in formulating a general purpose Helmholtz solver in two dimensions. The interpretation of the term `stabilized' here is twofold: (1) the condition number of the interpolation matrix is significantly reduced and (2) the system matrix in the linear system has stabilized eigenvalue spectra irrespective of any irregularity in the stencil, or the size of the stencil. The use of this hybrid RBF circumvents the ill-conditioning in RBF-FD and provides stabilized evaluation at a reduced computational cost, which is practically equal to that obtained by RBF-Direct methods. Although we have selected the kernel parameters based on the numerical tests in the work, future work involving finer tuning of the parameters is likely to further improve the presented algorithm. The absorbing boundary condition works excellent with the proposed RBF-FD discretization as we observe no significant spurious reflections from the truncated computational boundary. However, further improvement can be made in the presented approach by incorporating perfectly matched layers (PML) to increase the accuracy. Also, efficient ways to determine the kernel-parameters in the hybridization can be added in the present approach.

\begin{figure}
\centering
\includegraphics[scale=0.55]{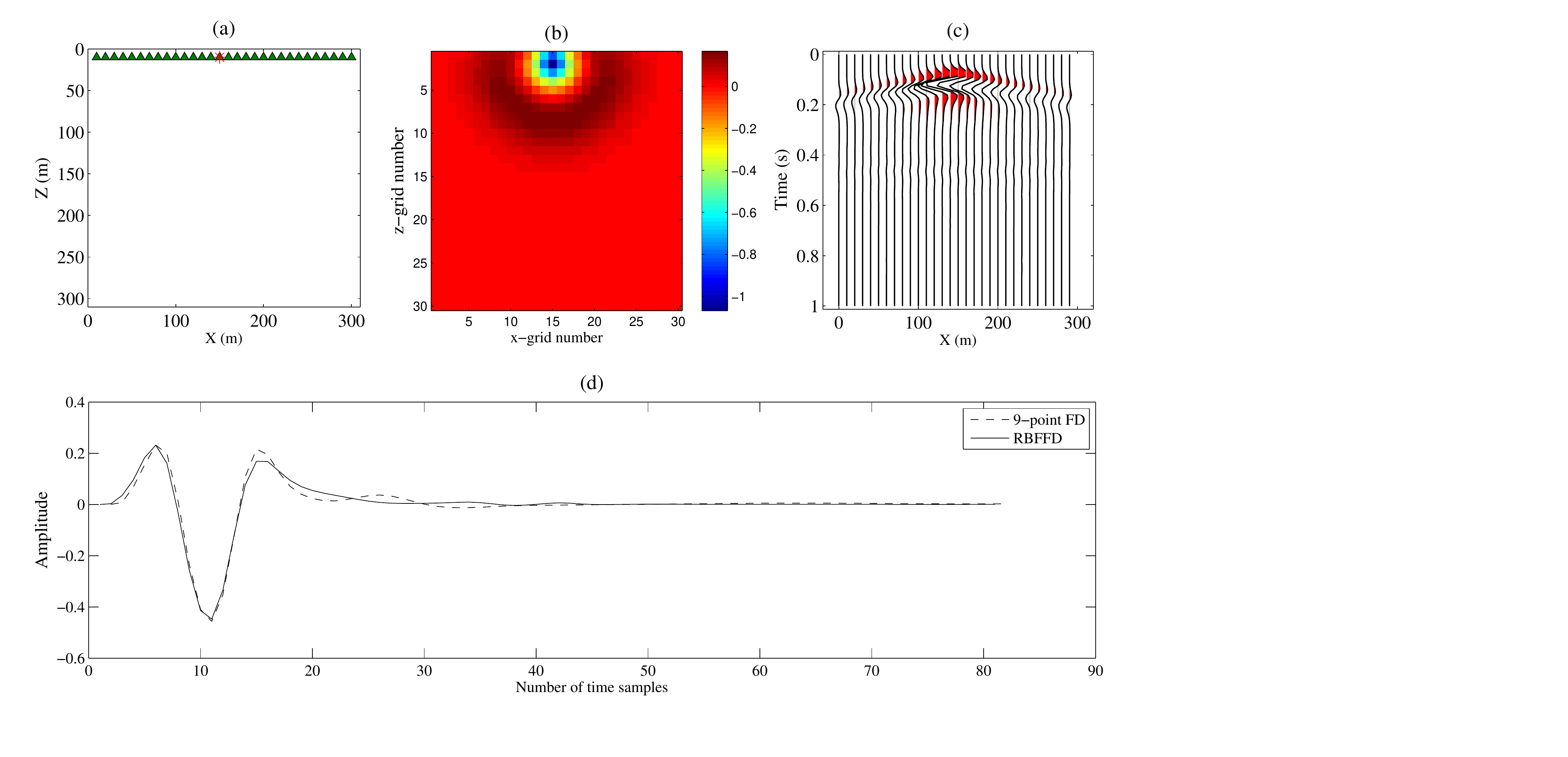}
\caption{ (a) The computational domain and the location of source and receivers, (b) acoustic wavefield at $t= 0.1750s$, (c) shot gather of seismograms, and (d) comparison of the seismogram passing through the source location, computed using RBF-FD and 9-point FD with 100m thick perfectly matched layers.}
\label{fig:seismo}
\end{figure}
\end{myfont}
\bibliographystyle{apalike}
{\bibliography{example}}
\label{lastpage}
\end{document}